\documentclass[12pt]{amsart}
\usepackage[francais]{babel}
\usepackage{graphics}
\usepackage{graphicx}
\usepackage{pdfpages}

\title{Des crit\`eres de  transcendance inspir\'es par un texte de Kolberg dat\'e de 1962} 

\author{Labib Haddad}
\address{120 rue de Charonne, 75011 Paris, France}
\email{labib.haddad@wanadoo.fr}

\usepackage{amssymb}
\usepackage{amsmath}

\usepackage[T1]{fontenc}

\newcommand{\su}{\subsection*}
\newcommand{\head}{\section*}
\newcommand{\noi}{\noindent}
\newcommand{\se}{\noi{\bf En effet}}

\newcommand{\Ž}{\'e}

\newcommand{\ˆ}{\`a}
\newcommand{\}{\`u}

\newcommand{\Q}{\mathbb Q}

\newcommand{\Z}{\mathbb Z}

\newcommand{\cal}{\mathcal}

\newcommand{\leqs}{\leqslant}

\newcommand{\geqs}{\geqslant}

\newcommand{\guil}{\guillemotleft}  
\newcommand {\guir}{\guillemotright}

\newcommand{\ali} {\begin{aligned}}   
\newcommand{\ala} {\end{aligned}}

\newcommand {\et}{\ \text{et}\ }

\newcommand {\ou}{\ \text{ou}\ }

\newcommand {\si}{\ \text{si}\ }

\newcommand {\pour}{\ \text{pour}\ }

\newcommand {\pourtout}{\ \text{pour tout}\ }

\newcommand {\dou}{\ \text{d'o\}\ }

\newcommand{\stm}{\smallsetminus}

\newcommand{\lopar}{\noi \{$\looparrowright$ \ }

\newcommand{\bc}{\begin{cases}}
\newcommand{\ec}{\end{cases}}

\begin{document}
\maketitle

\thispagestyle{empty}

\markboth{Labib Haddad}{Oddmund Kolberg}

Dans [0, 1962], {\sc  O. Kolberg} \Žnonce et \Žtablit le r\Žsultat suivant.

\su{Th\Žor\me} On se donne un entier $a>0$, un nombre rationnel $r$, distinct de $-1,-2,\dots, -a+1$, et un polyn\™me $P(z)$ non nul, \ˆ coefficients alg\Žbriques.  Alors, pour  $x$ alg\Žbrique tel que $0 < |x| < 1/e$, la somme de la s\Žrie suivante, $S$, est un nombre transcendant :
$$S = \sum^\infty_{n=1}\frac{(n+r)^{n-a}P(n)}{n!} x^n.$$
\su{Br\ve esquisse de la d\Žmonstration de Kolberg}  Pour chaque entier $k \in \Z$, on introduit la fonction :
$$f_k(x,y) = \sum^\infty_{n=0} \frac{x^{y+n}(y+n)^{n-k}}{n!}.$$
Nous poserons
$$T_k(x,y) = \sum^\infty_{n=0} (y+n)^{n-k}\frac{x^n}{n!},$$
de sorte que 
$$f_k(x,y) = x^y T_k(x,y).$$
Pour $r\neq 0$, en utilisant ces fonctions, et le changement de variable 
$$x = te^{-t}\ , \ |t|\leqs 1,$$ 
on montre que l'on a 
$$x^rS = t^rg(t)$$
o\ g(s) est une fraction rationnelle en $s$, dont les coefficients sont des nombres alg\Žbriques. D'apr\s le th\Žor\me de Lindemann, si $t$ est alg\Žbrique non nul, $x$ est transcendant. Comme $x$ est alg\Žbrique, le nombre $t$ est transcendant. On montre enfin que si $t$ est transcendant, le nombre $t^rg(t)$ l'est aussi. 

\noi Pour $r=0$, on a une d\Žmonstration semblable. \qed

\

La revue savante dans laquelle l'article de Kolberg est paru ne semble pas tr\s accessible [peu de num\Žros parus et cessation d'indexation apr\s 1972 dans {\bf zbMath}]. Aussi, donnerons-nous, ci-dessous,  dans un des derniers paragraphes, davantage de d\Žtails sur cette d\Žmonstration, en essayant d'en \Žclaircir certains points et de la simplifier.  

Notre propos, ici,  est  d'expliciter d'abord  sous une forme g\Žn\Žrale, un crit\re de transcendance que Kolberg utilise dans sa d\Žmonstration. D'introduire et d'\Žtudier  ensuite, sous le nom de {\bf quatuor}, des suites de fonctions qui g\Žn\Žralisent celles particuli\res des $f_k(x,y)$ et $T_k(x,y)$  qu'utilise Kolberg  dans sa d\Žmonstration, afin d'en \Žlargir le champ.
\

\head{Lemme (un crit\re de transcendance)}

\ 

Soit $g(s)$ une fraction rationnelle en $s$  dont les coefficients sont des nombres alg\Žbriques. On d\Žsigne par $E$ l'ensemble (exceptionnel) des nombres entiers $n\in \Z$ pour lesquels la fraction rationnelle $s^ng(s)$ est constante. On se donne un nombre rationnel $r \in \Q\stm E$. Alors, si  $t$ est un nombre transcendant, le nombre $t^rg(r)$ l'est  \Žgalement. 

\su{D\Žmonstration} Soient
$$g(s) =\frac{u(s)}{v(s)}\ , \ r= \frac{p}{q}$$
les formes irr\Žductibles de $g(s)$ et de $r$. Soit
$$t^r\frac{u(t)}{v(t)} = d.$$
On a
$$t^{p/q} u(t) = dv(t)$$
$$t^pu^q(t) = d^qv^q(t)$$
Ainsi, suivant que l'on a $r  >0, r=0, \ou r<0$, le nombre $t$ est racine du {\bf polyn\™me} :
$$s^pu^q(s) - d^qv^q(s) \ , \ u(s) - dv(s) \ , \ou u^q(s) - d^qs^{-p}v^q(s).$$
Aucun de ces polyn\™mes n'est nul car $r$ n'appartient pas \ˆ l'ensemble exceptionnel $E$.  Si $d$ \Žtait alg\Žbrique, ces polyn\™mes auraient des coefficients alg\Žbriques, et $t$ le serait aussi !\qed

\

\head{Quelques identit\Žs remarquables}

\

\noi Pour chaque entier $s\geqs 0$, on introduit le polyn\™me en $x$ suivant dont le degr\Ž est $\leqs s$  : 
 $$\boxed {P(x,s) = \sum_{n=0}^{s+1} (-1)^n\binom {s+1}{n}  (x+n)^s}$$

\noi On a, en particulier, 

$P(x,0) = 1-1  = 0$

$P(x,1) = x - 2(x+1) + (x+2) = 0$

$P(x,2) = x^2 - 3(x+1)^2 + 3(x+2)^2 - (x+3)^2=0$

$P(x,3) = x^3 - 4(x+1)^3 + 6(x+2)^3 - 4(x+3)^3  + (x+4)^3 = \ ?$

\

\noi Plus g\Žn\Žralement, des calculs simples conduisent aux identit\Žs  remarquables suivantes dont on aura \ˆ se servir dans la suite :
$$\boxed{P(x,s) = 0, \pourtout  s\geqs 0}.$$

\

\se,  le coefficient de $x^{s-q}$ dans le polyn\™me $(x+n)^s$ est \Žgal  \ˆ
$$\binom s q n^q$$
Le  coefficient $a(q,s)$ de $x^{s-q}$ dans le polyn\™me $P(x,s)$ est donc
$$a(q,s) = \sum_{n=0}^{s+1} (-1)^n\binom {s+1}{n}\binom s q n^q= \binom s q\sum_{n=0}^{s+1} (-1)^n\binom {s+1}{n} n^q.$$
On pose
$$b(q,s) = \sum_{n=0}^{s+1} (-1)^n\binom {s+1}{n} n^q \ \dou \ a(q,s) = \binom s q b(q,s).$$
En particulier, 
$$a(0,s) = b(0,s)= \sum_{n=0}^{s+1} (-1)^n\binom {s+1}{n}= (1-1)^{s+1} = 0.$$
Pour $0 < q \leqs s$, on \Žcrit 
$$b(q,s) = \sum_{n=0}^{s+1} (-1)^n\binom {s+1}{n} n^q = (s+1)\sum_{n=1}^{s+1} (-1)^n\binom {s}{n-1} n^{q-1}=$$
$$= -(s+1)\sum_{m=0}^{s} (-1)^m\binom {s}{m}m^{q-1}$$
autrement dit,
$$b(q,s)= -(s+1)b(q-1,s-1).$$
Or, $b(0,t) = 0 \pourtout t\geqs0$. Par r\Žcurrence, on  a $b(q,s) = 0$.
On a ainsi $a(q,s) = 0, \pour \ q\geqs 0 \et s\geqs 0$. Tous les coefficients du polyn\™me $P(x,s)$ sont nuls !\qed

\
 
Y a-t-il une mani\re plus simple d'\Žtablir ces identit\Žs remarquables ? 
C'est  probable, mais je ne le sais pas !

\

\head{S\Žrie associ\Že}

\
 
On se donne une s\Žrie enti\re en $x$
$$H(x) = \sum_{n=0}^\infty u_n\frac{x^n}{n!}.$$
On d\Žfinit la {\bf s\Žrie associ\Že}
$$G(t) = H(te^{-t}) = \sum_{n=0}^\infty u_n\frac{t^ne^{-nt}}{n!}= \sum_{n=0}^\infty v_k\frac{t^k}{k!}.$$
Cette derni\re s\Žrie est le d\Žveloppement de Taylor de la fonction $G(t)$ de sorte que  
$$v_k = \left[\frac{d^k}{dx^k}H(t)\right]{t=0}.$$ 
 En posant $E^k = \frac{d^k}{dt^k}$, il vient
$$c_{k,n} = \left[E^k \frac{t^ne^{-nt}}{n!}\right]_{t= 0}, \ \text{et on aura} \
v_k = \sum_{n=0}^\infty c_{k,n}u_n.$$
\su{Calculs} Successivement, on  a
$$\left[E^k\frac{t^n}{n!}\right]_{t=0}  =
{\bc 1 \si k=n
\\ 0 \si k\neq n\ec}$$
$$E^k e^{-nt} = (-1)^k n^k e^{-nt}$$
$$[E^k e^{-nt}]_{t=0} =  (- n)^k , \text{avec la convention} \ 0^0  = 1$$
$$E^k \frac{t^n}{n!}e^{-nt} = \sum_{h=0}^k \binom  k h E^h\frac{t^n}{n!}E^{k-h}  e^{-nt}, \text{par la formule de Leibniz}$$
$$c_{k,n} = \left[E^k \frac{t^n}{n!}e^{-nt}\right]_{t=0}  = {\bc 
0 \si k < n\\
\binom k n (-n)^{k-n} \si n\leqs k\ec}$$
$$v_k  =  \sum_{n=0}^k \binom k n (-n)^{k-n} u_n.$$ 
Ainsi, en particulier, on a $v_0 = u_0$. Pour $k\geqs 1$, on a $\binom k 0 = 0$, donc
$$\boxed{v_k  =  \sum_{n=1}^k \binom k n (-n)^{k-n} u_n}$$
autrement dit 
$$\boxed{v_k = u_k  + \dots+  \binom k n (-n)^{k-n} u_n + \dots  +    (-1)^{k-1}ku_1.}$$
On v\Žrifie alors que  les  expressions suivantes donnent, r\Žciproquement, les $u_n$ en fonction des $v_h$.
$$\boxed{u_n = \sum_{h=1}^n \binom {n -1}{ h-1} n^{n-h}v_h}$$
autrement dit
$$\boxed{u_n = v_n +  \dots + \binom {n-1}{h-1} n^{n-h} v_h + \dots+ n^{n-1}v_1.}$$

\

\su{V\Žrification}  Il s'agit de v\Žrifier que l'on a
$$v_k = \sum_{n=1}^k \binom k n (-n)^{k-n}\sum_{h=1}^n \binom {n-1} {h -1} n^{n-h}v_h.$$
C'est une combinaison lin\Žaire de $v_1,v_2,\dots,v_n,$
$$d_1v_1 + \dots + d_hv_h + \dots d_nv_n$$
o\  le coefficient de $v_h$ est
$$d_h = \sum_{n=h}^k \binom k n\binom {n-1}{h-1}(-n)^{k-n}n^{n-h}.$$ 
On a $d_k = 1$. 

\noi  Pour $1 \leqs h < k$, on v\Žrifie que $d_h  = 0$ en se servant des identit\Žs remarquables, $P(x,s) = 0$, pr\Žsent\Žes ci-dessus. En effet, on a 
$$d_h = (-1)^k\sum_{n=h}^k(-1)n\frac{k!(n-1)!}{n!(k-n)!(h-1)!(n-h)!}n^{k-h}$$
$$=(-1)^k\frac{k!}{(h-1)!(k-h)!}\sum_{n=h}^k (-1)^n\binom {k-h}{n-h}n^{k-h-1}.$$
On pose $k-h = s+1$. Il vient
$$\sum_{n=h}^k (-1)^n\binom {k-h}{n-h}n^{k-h-1}= \sum_{n=h}^{h+s+1}(-1)^n\binom {s+1}{n-h}n^s =$$
 $$= (-1)^h\sum_{n=0}^{s+1}(-1)^n\binom {s+1} n(h+n)^s = (-1)^hP(h,s) = 0.$$\qed

\su{Insistons en le disant encore autrement} Pour que
$$G(t) = \sum_{n=0}^\infty v_n\frac{t^n}{n!},$$
soit la s\Žrie associ\Že \ˆ 
$$H(x) = \sum_{n=0}^\infty u_n\frac{x^n}{n!},$$
il faut  et il suffit que l'on ait $u_o = v_0$ et, pour  $n\geqs 1$,
$$\boxed{u_n = v_n +  \dots + \binom {n-1}{m-1} n^{n-m} v_m + \dots+ n^{n-1}v_1,}$$
ce qui \Žquivaut encore \ˆ :
$$\boxed{v_n = u_n  + \dots+  \binom n m (-n)^{n-m} u_m + \dots  +    (-1)^{n-1}nu_1.}$$

\

\su{Un exemple prototype} La s\Žrie associ\Že \ˆ  
$$H(x) = \sum_{n=0}^\infty u_n\frac{x^n}{n!}, \ u_n =(y+n)^{n-1},$$
est la s\Žrie 
$$G(t) = \sum_{n=0}^\infty v_n\frac{t^n}{n!}, \ v_n = y^{n-1}, \ \text{autrement dit,}$$

$$G(t) =\frac{e^{yt}}{y}.$$

\se, on a $u_0 = v_0 = 1/y$ et
$$u_n = \sum_{h=1}^n \binom {n -1}{ h-1} n^{n-h}y^{h-1} = (y+n)^{n-1}.$$
\qed

\

On convient d'appeler {\bf suite} toute famille ind\Žx\Že par l'ensemble $\Z$ des entiers.

\

\head{Les quatuors}

\

\su{D\Žfinition} Un {\bf quatuor} est compos\Ž de quatre suites, $F, G, H, K,$ o\ 
$$K=(K_k(x,y))_{k\in \Z}, H  = (H_k(x,y))_{k\in \Z},$$
$$F = (F_k(t,y))_{k\in\Z}, G= (G_k(t,y))_{k\in\Z}.$$
Ce sont quatre suites de fonctions. Ces fonctions doivent \tre li\Žes par les relations suivantes
$$K_k(x,y) = x^yH_k(x,y)$$
$$G_k(t,y) = H_k(te^{-t},y)$$
$$F_k(t,y) = K_k(te^{-t},y)$$ 
et satisfaire les condition suivantes :
\[x\frac{d}{dx}K_{k+1}(x,y) = K_k(x,y)\tag*{(1k)}\]

\noi Il s'ensuit, comme on le v\Žrifie par des calculs simples, que l'on aura alors l'ensemble de  toutes les relations que voici  :

\

  $\bc

K_k(x,y) = x^yH_k(x,y)\\ 
\\

G_k(t,y) = H_k(te^{-t}, y)\\ 
\\

F_k(t,y) = K_k(te^{-t},y) =  t^ye^{-yt}G_k(t,y)\\
\\

\displaystyle
(1k) \ x\frac{d}{dx}K_{k+1}(x,y) = K_k(x,y)\\
\\

\displaystyle
(2k) \ K_{k+1}(x,y) = \int_0^x \frac{K_k(z,y)}{z}dz\\

\\
\displaystyle
(3k) \ x\frac{d}{dx} H_{k+1}(x,y) + yH_{k+1}(x,y) = H_k(x,y)\\
\\

\displaystyle
(4k) \ x^yH_{k+1}(x,y) =  \int_0^x  z^{y-1}H_k(z,y)dz \\
\\

\displaystyle
(5k) \ \frac{t}{1-t}\frac{d}{dt} G_{k+1}(t,y) + yG_{k+1}(t,y) = G_k(t,y)\\
\\

\displaystyle
(6k) \ t^ye^{-yt}G_{k+1}(t,y) = \int_0^t \frac{1-z}{z}z^ye^{-yz}G_k(z,y)dz\\
\\

\displaystyle
(7k) \ \frac{d}{dt}F_{k+1}(t,y) = \frac{1-t}{t}F_k(t,y)\\
\\

\displaystyle
(8k) \ F_{k+1}(t,y) = \int_0^t \frac{1-z}{z} F_k(z,y) dz\\
\\
 
\ec$ 

\

\

\

\noi En particulier, les relations (3k) sont \Žquivalentes aux relations (1k). On s'en souviendra, \ˆ loccasion.

\su{Pour ce qui est des d\Žrivations !} En outre, on suppose qu'il existe un intervalle non vide, $V = [0, \varepsilon[$, voisinage de $0$ \ˆ droite, tel que chacune des fonctions $H_k(x,y)$ poss\de une d\Žriv\Že par rapport \ˆ $x$ en tout  point de $V$. Il s'ensuit que toutes les fonctions, $H_k(x,y) \et K_k(x,y)$, sont ind\Žfiniment d\Žrivables par rapport \ˆ $x$ en tout point de $V$. De m\me, toutes les fonctions $G_k(t,y) \et F_k(t,y)$ sont ind\Žfiniment d\Žrivables par rapport \ˆ $t$ en tout point d'un voisinage de $0$ \ˆ droite, $W = [0,\eta[$.

\su{Unicit\Ž et conditions initiales} En partant de la fonction $K_0(x,y)$ toute seule, on  peut reconstituer le quatuor tout entier ! \se, on obtient $K_k(x,y)$, pour $k = 1,2,3,\dots$, par int\Žgrations successives :
$$K_{k+1}(x,y) = \int_0^x \frac{K_k(z,y)}{z}dz.$$
On obtient $K_k(x,y)$, pour $k = -1, - 2,-3,\dots$, par d\Žrivations successives :
$$K_k(x,y) = x\frac{d}{dx}K_{k+1}(x,y).$$
Les fonctions $H_k, G_k \et F_k$ s'ensuivent.\qed

\

\noi Il en irait de m\me si l'on prenait la seule fonction $K_h(x,y)$, par exemple, au lieu de $K_0(x,y)$. Un simple d\Žcalage ! Dans ce sens, on dira que la  fonction $K_h(x,y)$,  \ˆ elle seule, {\bf engendre} le quatuor. Chacune des fonctions $K_k(x,y)$ est ainsi un {\bf g\Žn\Žrateur} du quatuor !
 
\

\su{L'exemple princeps. Le quatuor  de Kolberg} 

\

\

\noi C'est le quatuor $\cal K = (F,G,H,K)$  o\ 
$$H_k(x,y) = T_k(x,y) = \sum^\infty_{n=0} (y+n)^{n-k}\frac{x^n}{n!}.$$
$$K_k(x,y) = x^yT_k(x,y) = f_k(x,y).$$
Les $T_k(x,y)$ sont des s\Žries enti\res en $x$ qui ont toutes le m\me rayon de convergence, $1/e$. En effet, en  posant $u_n = (y+n)^{n-k}/n!$, on v\Žrifie que la limite de $u_{n+1}/u_n$ est \Žgale \ˆ $e$.

\

\su{Permanence} On peut dire, en bref :  toute {\bf combinaison lin\Žaire} de quatuors  est un quatuor. En particulier {\bf la diff\Žrence} de deux quatuors est un quatuor. 
De m\me, les d\Žcal\Žs d'un quatuor sont des quatuors.

\su{Explicitons, s'il en est besoin} 
\

\noi Soit une famille donn\Že de quatuors $\cal Q^j= $ 
$$K^j=(K^j_k(x,y))_{k\in \Z}, H^j  = (H^j_k(x,y))_{k\in \Z},$$
$$F^j = (F^j_k(t,y))_{k\in\Z}, G^j= (G^j_k(t,y))_{k\in\Z}.$$
La combinaison lin\Žaire $\lambda_1\cal Q^1 +\dots \lambda_p\cal Q^p$ est le quatuor $\cal Q=$
$$(K_k(x,y) = \lambda_1K^1_k(x,y)+ \dots + \lambda_pK^p(x,y))_{k\in Z}$$
$$(H_k(x,y) = \lambda_1H^1_k(x,y)+ \dots + \lambda_pH^p(x,y))_{k\in Z}$$
$$(F_k(t,y) = \lambda_1F^1_k(t,y)+ \dots + \lambda_pF^p(t,y))_{k\in Z}$$
$$(G_k(x,y) = \lambda_1G^1_k(t,y)+ \dots + \lambda_pG^p(t,y))_{k\in Z}$$
La diff\Žrence $\cal Q^1 - \cal Q^2$ est le quatuor $\cal Q =$
$$(K^1_k(x,y) - K^2_k(x,y))_{k\in \Z}, (H^1_k(x,y) - H^2_k(x,y))_{k\in \Z}$$
$$F^1_k(t,y) - F^2_k(t,y))_{k\in \Z}, (G^1_k(t,y) - G^2_k(t,y))_{k\in \Z}.$$
Le d\Žcal\Ž d'odre $d$ du quatuor  $\cal Q$ est le quatuor $\cal Q^{\to d}=$
$$K^{\to d}=(K_{k+d}(x,y))_{k\in \Z}, H^{\to d}  = (H_{k+d}(x,y))_{k\in \Z},$$
$$F^{\to d} = (F_{k+d}(t,y))_{k\in\Z}, G^{\to d}= (G_{k+d}(t,y))_{k\in\Z}.$$

\

\head{Une classe tr\s particuli\re de quatuors}

\

On va explorer la classe sp\Žciale des quatuors $\cal Q = (F, G, H, K)$ o\ les $H_k$ et les $G_k$ sont des s\Žries enti\res :
$$H_k(x,y) = \sum_{n=0}^\infty u_{k,n}\frac{x^n}{n!}, \  G_k(t,y) = \sum_{n=0}^\infty v_{k,n}\frac{t^n}{n!},$$
et les coefficients, $u_{k,n} \et v_{k,n}$, sont des fonctions de $y$.

\noi Pour qu'il en soit ainsi, il suffit que l'une des fonctions $H_k$ ou $G_k$ soit une s\Žrie enti\re car alors elles le seront toutes. La s\Žrie $G_k$ est l'associ\Že de la s\Žrie $H_k$ (voir ci-dessus) on a donc les relations suivantes : $u_{k,0} = v_{k,0}$ et, pour $n\geqs 1$ :
$$\boxed{u_{k,n} = v_{k,n} +  \dots + \binom {n-1}{m-1} n^{n-m} v_{k,m} + \dots+ n^{n-1}v_{k,1}}$$
$$\boxed{v_{k,n} = u_{k,n}  + \dots+  \binom n m (-n)^{n-m} u_{k,m} + \dots  +    (-1)^{n-1}nu_{k,1}.}$$

\

\su{La diff\Žrence avec un  quatuor quelconque.} Dans un quatuor quelconque,  $\cal Q = (F,G,H,K)$, chacune des fonctions  $H_k$ et $G_k$,  est ind\Žfiniment d\Žrivable sur un voisinage \ˆ droite $V$ de $0$ ; il existe donc bien une s\Žrie de Taylor, au point $0$, attach\Že \ˆ cette fonction. Cependant, le rayon de convergence de la s\Žrie peut \tre nul, et m\me s'il ne l'\Žtait pas, la somme de la s\Žrie pourrait ne pas \tre \Žgale \ˆ la fonction !

\

\noi {\bf Pour ces quatuors sp\Žciaux, la fonction est \Žgale \ˆ la somme de sa s\Žrie de Taylor !} C'est ce qui les caract\Žrise !

\

\noi {\it Cum grano salis }! On dira que ce sont des {\bf quatuors s\Žriels} !

\

\noi Bien entendu, chaque fonction $K_k(x,y)$ est somme de termes  en $x^yx^n$; c'est une s\Žrie mais ce n'est donc pas une s\Žrie enti\re sauf si $y$ est entier. De m\me, chaque fonction $F_k(t,y)$ est somme de  termes en $t^ye^{-ty}t^n$.

\

\

\head{Le quatuor de Kolberg}

\

C'est le quatuor s\Žriel $\cal K = (F,G,H,K)$  o\ 
$$H_k(x,y) = T_k(x,y) = \sum^\infty_{n=0} (y+n)^{n-k}\frac{x^n}{n!}.$$ 
$$K_k(x,y) = x^yT_k(x,y) = f_k(x,y).$$
Notons tout de suite ceci  : on a $T_k(0,y) = y^{-k}$, pour chaque $k\in\Z$. Tous les coefficients de ces s\Žries sont des fonctions {\bf analytiques} de la variable $y$. Les fonctions, $F_k(t,y), G_k(t,y), K_k(x,y), H_k(x,y)$, (limites uniformes de fonctions analytiques) sont donc toutes des fonctions analytiques en $y$.  Attention,  lorsque $k\geqs1$, ces fonctions ont un p\™le au point $y=0$. 

\

\noi On observe que l'on a bien
$$x\frac{d}{dx} T_{k+1}(x,y) + yT_{k+1}(x,y) = T_k(x,y).$$
En effet, en comparant les coefficients  des termes en $x^n/n!$ de ces trois s\Žries, on a bien $n(y+n)^{n-1-k} + y(y+n)^{n-k-1} = (y+n)^{n-k}$.

\noi On le sait d\Žj\ˆ, cela \Žquivaut \ˆ
$$x\frac{d}{dx} f_{k+1}(x,y) = f_k(x,y).$$

\noi Comme on l'a vu ci-dessus dans l'exemple prototype de s\Žrie associ\Že, on a 
$$G_1(t,y) =\frac{e^{yt}}{y} \ \dou F_1(t,y) = \frac{t^y}{y}.$$
Alias
$$T_1(te^{-t},y) = \frac{e^{yt}}{y} \ \et f_1(te^{-t},y) = \frac{t^y}{y}.$$
Ici, les sommes des s\Žries sont sous formes closes. Cela simplifie les calculs des d\Žriv\Žes et des primitives !

$$\boxed{\ \text{Pour} \ y\neq 0}$$

\su{La suite des $\mathbf{F_k(t,y)}$} On obtient toutes les $F_k(t,y)$ \ˆ partir de $F_1(t,y) = \frac{t^y}{y}$ en utilisant les formules (7k) et (8k) :
$$F_k(t,y) = \frac{t}{1-t}\frac{d}{dt}F_{k+1}(t,y) , \ F_{k+1}(t,y) = \int_0^t \frac{1-z}{z} F_k(z,y) dz.$$
Pour $k = 0, -1, -2, \dots$, on obtient $F_k(x,y)$ par d\Žrivations successives et, pour $k = 2, 3, 4, \dots$, par int\Žgrations successives.
Ainsi, par exemple,
$$F_0(t,y) = \frac{t^y}{1-t}$$
$$F_{-1}(t,y) = \frac{t^y}{(1-t)^2}\left(y-1 + \frac{1}{1-t} \right)$$
 Plus g\Žn\Žralement, pour $k = 0, -1, -2, \dots$, par r\Žcurrence, on a
 $$F_k(t,y) = t^y(1-t)^{k-1}P_k(t,y)$$
 o\ $P_k(t,y)$ est un polyn\™me en $\frac{1}{1-t}$ de degr\Ž $-k$ dont les coefficients sont des polyn\™mes en $y$ \ˆ coefficicents rationnels :
 $$P_k(t,y) \in \Q[y]\left [\frac{1}{1-t}\right].$$
 De m\me, par exemple,
$$F_2(t,y) = \frac{t^y}{y}\left(\frac{1}{y} - \frac{t}{y+1}\right).$$
$$F_3(t,y) = \frac{t^y}{y}\left\{\frac{1}{y^2} - \left(\frac{1}{y}+\frac{1}{y+1}\right) \frac{t}{y+1} + \frac{t^2}{(y+1)(y+2)}\right\}$$
et pour $k =1,2,3,\dots$, et $y\neq  -1, -2, -3,\dots, -k+1$, on a
$$F_k(t,y) = t^yQ_k(t,y)$$
o\ $Q_k(t,y)$ est un polyn\™me en $t$ de degr\Ž $k-1$ dont les coefficients sont des fractions rationnelles en $y$ \ˆ coefficients rationnels :
$$Q_k(t,y) \in \Q(y)[t]$$
Pour ($k \leqs 0$) et  pour ($k\geqs 1$ mais $y\neq  -1, -2,\dots, -k +1$), on a ainsi :
$$(9) \ \ \quad \boxed{\ F_k(t,y) =t^yR_k(t,y) \ \et R_k(t,y) \in \Q(y)(t)}$$
autrement dit,  $R_k(t,y)$ est une fraction rationnelle en $t$ dont les coefficients sont des fractions rationnelles en $y$.

\

{\bf  Tout cela est bien mentionn\Ž par Kolberg. Rien de nouveau.}

$$\boxed{\ \text{Pour} \ y = 0}$$

\

Afin d'\Žviter le p\™le en $y=0$, on introduit un nouveau quatuor $\cal K^\vee$ enngendr\Ž par la s\Žrie
$$K_1^\vee(x) =\sum^\infty_{n=1} n^{n-1}\frac{x^n}{n!}.$$
On va accorder un peu d'attention \ˆ ce quatuor s\Žriel, sorte de quatuor opus 2 de Kolberg, qui ne d\Žpend pas de la variable $y$ et dans lequel on a donc
$$K_k^\vee(x) = H_k^\vee(x).$$
Partant de $K_1^\vee(x)$, les formules (1k) et (2k), donnent :
$$K_k^\vee(x) = H_k^\vee(x) = \sum^\infty_{n=1} n^{n-k}\frac{x^n}{n!}, \pourtout k\in \Z.$$
D'autre part, soit 
$$G_1^\vee(t) = \sum_{n=1}^\infty v_n\frac{t^n}{n!}$$
la s\Žrie associ\Že \ˆ $H_1^\vee(x)$.
Comme dans l'exemple prototype,  si l'on prend  $u_n = n^{n-1}$, on v\Žrifie que $v_1 = 1$  et $v_n = 0, \pourtout n>1$. Ainsi, 
$$G_1^\vee(t) = t \dou F_1^\vee(t) = t.$$
\`A l'aide des formules (7k) et (8k), comme dans le cas g\Žn\Žral, on obtient la suite de toutes les  $F_k^\vee(t)$ \ˆ partir de $F_1^\vee(t) = t$.
$$F_k^\vee(t) = \frac{t}{1-t}\frac{d}{dt}F_{k+1}^\vee(t) , \ F_{k+1}^\vee(t) = \int_0^t \frac{1-z}{z} F_k^\vee(z) dz.$$
On constate que, pour $k \leqs 0$, la fonction $F_k^\vee(t)$ est une fraction rationnelle en $t$ \ˆ coefficients rationnels; pour  $k >0$, c'est un polyn\™me en $t$ \ˆ coefficients rationnels. Dans les deux cas, on aura 
$$(10) \qquad \boxed{F_k^\vee(t) \in \Q(t).}$$

\

\head{D\Žmonstration du  Th\Žor\me de  Kolberg}

\

\su{Un rappel} Soit $x = te^{-t}$. D'apr\s Lindemann, si $t$ est alg\Žbrique non nul, $x$ est transcendant. Donc, si $x$ est alg\Žrique non nul, $t$ est transcendant !

\

\noi On se donne un entier $a>0$, un polyn\™me non nul, $P(z)$, \ˆ coefficients   alg\Žbriques, et un nombre rationnel $r \neq -1, -2, \dots, -a+1$. On a 
$$S = \sum^\infty_{n=1}\frac{(n+r)^{n-a}P(n)}{n!} x^n.$$
En posant
$$u_n = \frac{(n+r)^{n-a}P(n)}{n!}$$
on voit (facilement) que la limite de $u_{n+1}/u_n$ est \Žgale \ˆ $e$ de sorte que le rayon de convergence de la s\Žrie $S$ est \Žgal \ˆ $1/e$.

$$\boxed{\ \text{Le cas} \ r\neq 0}$$

\

Il est loisible de mettre le polyn\™me $P(z)$ sous la forme :
$$P(z) = \sum_{j=q}^mp_j(z+r)^j$$
o\ les coefficients $p_j$ sont des nombres alg\Žbriques, $q\leqs m$, $p_q \neq 0$, $p_m\neq 0$, et $m$ est le degr\Ž de $P(z)$.  On a ainsi :
$$S = \sum_{n=1}^\infty\sum_{j=q}^mp_j(n+r)^{n-a+j}\frac{x^n}{n!}.$$
En faisant le changement de variable $k = a -j$, $j = a- k$, il vient :
$$S = \sum_{n=1}^\infty\sum_{k = a-q}^{a-m}p_{a-k}(n+r)^{n-k}\frac{x^n}{n!}.$$
On pose alors $b = a -m \ , \ c=a-q \ , \ A_k = p_{a-k}$. Il vient :
$$S = \sum_{n=1}^\infty\sum_{k = b}^{c}A_k(n+r)^{n-k}\frac{x^n}{n!}$$ 
o\ $b\leqs c$, et $A_b \neq 0$, $A_c\neq 0$.

\noi On utilise les fonctions du quatuor de Kolberg :
$$f_k(x,y) = \sum^\infty_{n=0} \frac{x^{y+n}(y+n)^{n-k}}{n!}.$$
Pour $b\leqs c$ des entiers quelconques dans $\Z$,  $A_k$ des nombres alg\Žbriques, et $A_b \neq 0$, $A_c \neq 0$, on introduit la combinaison lin\Žaire suivante :
$$L = \sum_{k=b}^cA_kf_k(x,r)$$
Il vient
$$ L = \sum_{n=0}^\infty\sum_{k=b}^c A_k(n+r)^{n-k}\frac{x^{n+r}}{n!}= x^r\sum_{n=0}^\infty\sum_{k=b}^c A_k(n+r)^{n-k}\frac{x^n}{n!} $$
Autrement dit, on a $L = x^rS$.

\noi Avec le changement de variable $x = te^{-t}$, par la relation (9) signal\Že ci-dessus, la combinaison lin\Žaire $L$ est de la forme :
$$L = t^rg(t)$$
o\ $g(t)$ est une fraction rationnelle en $t$ dont les coefficients sont des nombres alg\Žbriques. 

\noi Il est clair que $t^rg(t)$ n'est pas contante comme fonction de $t$. D'apr\s le crit\re de transcendance \Žnonc\Ž comme lemme ci-dessus, si $t$ est transcendant, alors $L = t^rg(t)$ est transcendant. Si le nombre non nul $x$ est alg\Žbrique, le nombre $t$ est transcendant et le nombre  $x^r$ est alg\Žbrique, de sorte que $S$ est transcendant puisque $L=x^rS$ l'est. \qed

$$\boxed{\ \text{Le cas} \ r = 0}$$

\

On proc\de de mani\re analogue, mais en utilisant l'opus 2. Le polyn\™me $P(z)$ est de forme :
$$P(z) = \sum_{j=q}^mp_jz^j$$
o\ les coefficients $p_j$ sont des nombres alg\Žbriques, $q\leqs m$, $p_q \neq 0$, $p_m\neq 0$, et $m$ est le degr\Ž de $P(z)$.  On a ainsi :
$$S = \sum_{n=1}^\infty\sum_{j=q}^mp_jn^{n-a+j}\frac{x^n}{n!}.$$
En faisant le changement de variable $k = a -j$, $j = a- k$, il vient :
$$S = \sum_{n=1}^\infty\sum_{k = a-q}^{a-m}p_{a-k}n^{n-k}\frac{x^n}{n!}.$$
On pose alors $b = a -m \ , \ c=a-q \ , \ A_k = p_{a-k}$. Il vient :
$$S = \sum_{n=1}^\infty\sum_{k = b}^{c}A_kn^{n-k}\frac{x^n}{n!}$$ 
o\ $b\leqs c$, et $A_b \neq 0$, $A_c\neq 0$.

\noi On utilise les fonctions du quatuor opus 2  :
$$f_k^\vee(x) = \sum^\infty_{n=1} \frac{n^{n-k}x^n}{n!}.$$
Pour $b\leqs c$ des entiers quelconques dans $\Z$,  $A_k$ des nombres alg\Žbriques, et $A_b \neq 0$, $A_c \neq 0$, on introduit la combinaison lin\Žaire suivante :
$$L = \sum_{k=b}^cA_kf_k^\vee(x)$$
Il vient
$$ L = \sum_{n=0}^\infty\sum_{k=b}^c A_kn^{n-k}\frac{x^n}{n!}= \sum_{n=0}^\infty\sum_{k=b}^c A_kn^{n-k}\frac{x^n}{n!} $$
Autrement dit, on a $L = S$.

\noi Avec le changement de variable $x = te^{-t}$, par la relation (10) signal\Že ci-dessus, la combinaison lin\Žaire $L$ est de la forme :
$$L = g(t)$$
o\ $g(t)$ est une fraction rationnelle en $t$ dont les coefficients sont des nombres alg\Žbriques. 

D'apr\s le crit\re de transcendance \Žnonc\Ž comme lemme ci-dessus, si $t$ est transcendant, alors $L = g(t)$ est transcendant donc $S$ est transcendant.

Si le nombre non nul $x$ est alg\Žbrique, le nombre $t$ est granscendant de sorte que $S$ est transcendant puisque $L=S$ l'est. \qed

\

\

\

\

\head{Remarques compl\Žmentaires}

\

Kolberg fait encore observer ceci.

\su{1} Le th\Žor\me est encore  vrai pour les entiers $a\leqs 0$. En effet, soit $a\leqs  0$. On \Žcrit :
$$(n+r)^{n-a} P(n) = (n+r)^{n-1}(n+r)^{-a +1} P(n)= (n+r)^{n-1}Q(n)$$
o\ $(-a+1) \geqs 0$ et $Q(n) =(n+r)^{-a +1} P(n)$ est un polyn\™me en $n$.

\su{2} Pour $k>0$, la s\Žrie
$$f_k^\vee(x) = \sum^\infty_{n=1} \frac{n^{n-k}x^n}{n!},$$
dont le rayon de convergence est $1/e$ est convergente pour $x =1/e$ qui correspond \ˆ $t=1$. Il s'ensuit que
$$\sum^\infty_{n=1} \frac{n^{n-k}e^{-n}}{n!}$$
est un nombre rationnel !

\

\head{{ \emph{\Large Kolbergisation}}}

\

Le proc\Žd\Ž utilis\Ž dans la d\Žmonstration du th\Žor\me, \ˆ deux reprises, peut s'appliquer \ˆ d'autres cas ! C'est un  proc\Žd\Ž qu'on devrait appeler 
\textbf{\textit {Kolbergisation}} en hommage \ˆ son auteur, bien entendu !

En voici un exemple, un peu au hasard. 

\noi Introduisons le quatuor $\cal E =(F,G,H,K)$ o\ $G_0(t) = t^2/2$  [un quatuor en $t$ mineur, pour ainsi dire]. Comme on le fait dans la d\Žmonstration pour le cas $r=0$, on reconstitue le quatuor. En prenant
$$H_0(x) = \sum_{n=0}^\infty u_n\frac{x^n}{n!}, \ G_0(t) = \sum_{n=0}^\infty v_n\frac{t^n}{n!},$$
il vient $v_2 = 1$ et, pour $n\neq 2$, $v_n = 0$, donc
$$u_n = \binom {n-1} 1 n^{n-2}v_2 = (n-1)n^{n-2}$$
$$H_0(x) = \sum_{n=1}^\infty (n-1)n^{n-2}\frac{x^n}{n!}$$
$$K_k(x) = H_k(x) = \sum_{n=1}^\infty (n-1)n^{n+k-2}\frac{x^n}{n!}$$
Le rayon de convergence des ces s\Žries est $1/e$. On v\Žrifie simplement que $F_k(t)$ est une fraction rationnelle en $t$ \ˆ coefficients rationnels (comme dans la d\Žmonstration ci-dessus pour le cas $r=0$) 

\

Soit alors
$$L = \sum_k A_k K_k(x)$$
une combinaison lin\Žaire quelconque des $K_k(x)$, \ˆ coefficients alg\Žbriques. Pour chaque nombre $x$ alg\Žbrique tel que $0 < |x| < 1/e$, le nombre $L$ est transcendant ! Un exemple de kolbergisation !

\

\noi Mais que sont ces combinaisons lin\Žaires $L$ ?

\

\noi Ce sont (sauf  erreur) les s\Žries enti\res de la forme :
$$\boxed{\sum_{n=2}^\infty (n-1)n^{n+a}  P\left(\frac{1}{n}\right) \frac{x^n}{n!}}$$
o\ $a$ est un entier quelconque et $P(z)$ un polyn\™me non  nul \ˆ coefficients alg\Žbriques. Un r\Žsultat que Kolberg aurait pu \Žtablir !

\

\head{En r\Žsum\Ž}

\

Dans un article, voir  [0],  dat\Ž de 1962, Kolberg \Žnonce et \Žtablit un th\Žor\me sur la transcendance des valeurs des sommes d'une classe de certaines s\Žries enti\res en $x$, pour les valeurs alg\Žbriques de $x$. Il s'appuie sur le th\Žor\me de Lindemann et utilise, en passant, un crit\re de transcendance pour les valeurs de certaines fractions rationnelles. On explicite ce dernier crit\re. On \Žclaircit certains points d\Žlicats de la d\Žmonstration et on montre comment on peut \Žtendre son champ. 

\

\head{Addendum}

\

Arriv\Ž \ˆ un certain point de sa d\Žmonstration, Kolberg \Žcrit cette phrase quelque peu sibylline :

\

\lq\lq To obtain an expression for $f_k(x,y)$ we use a well known formula due to Legendre (see [1]). In fact, putting
$$x = te^{-t} , |t| \leqs 1$$
we have
$$f_1(x,y) = \frac{x^y e^{yt}}{y} = \frac{t^y}{y}."$$
Le [1] est un renvoi \ˆ A. M. Legendre, {\sl Exercices de calcul int\Žgral}, Paris (1811), sans aucune indication de la page o\ se trouve cette formule. Le livre de Legendre est compos\Ž de plusieurs tomes dont le premier, de 386 pages,  comporte  trois parties : Des fonctions elliptiques. Des int\Žgrales eul\Žriennes. Des quadratures. Sans compter le tome  intitul\Ž Quatri\me partie, en 152 pages, et le Suppl\Žment \ˆ la premi\re partie, en 50  pages. Autant chercher une aiguille dans une botte de foin !

\

\lopar Cela fait penser, immanquablement, \ˆ la conf\Žrence faite par Jean-Pierre Serre, en 2003, au Harvard's Basic Notions seminar, intitul\Že How to write mathematics badly [2]. Il en existe une transcription \Žcrite par Maxine Calle [3].

Serre y moque malicieusement ce travers, parmi plusieurs mauvaises mani\res d'\Žcire des math\Žmatiques : renvoyer \ˆ un livre ou manuel de plusieurs centaines de pages, sans donner aucune indication de la page ou du num\Žro du th\Žor\me auquel on se r\Žf\re !\}

\

Je ne sais toujours pas de quelle formule \guil bien connue\guir\ de Legendre il s'agit  ! J'ai pourtant beaucoup cherch\Ž, sans aucun succ\s, en vain. 

\
 
En y r\Žfl\Žchissant un peu, on peut imaginer que Kolberg a simplement voulu, apr\s tout, user d'une mystification !

\

\head{Bibliographie}

\

\noi [0] {\sc O. Kolberg}, {\sl A class of power series with transcendental sums for algebraic values of the variable}. \AA rbok Univ. Bergen, Mat.-Naturv. Ser. (1962) No. 18, 6 p.

\  

\noi  [1]  {\sc A.  M. Legendre}, {\sl Exercices de calcul int\Žgral},   Mme veuve COURCIER,  Imprimeur-Libraire, quai des Augustins, Paris (1811).

\

\noi [2] {\sc J.-P. Serre}, {\sl How to write mathematics badly}. 

\noi Video: {\tt https://www.youtube.com/ watch?v=ECQyFzzBHlo}.

\

\noi [3] {\sl How to write mathematics badly} by John Paul Serre, Transcribed by Maxine Cale {\tt http://web.sas.upenn.edu/callem/}

\

\

\

\enddocument